\magnification=\magstep1

\raggedright
\font\large=cmbx10 scaled \magstep1
\pageno=1
\normalbaselines 
\def\Const{{\it Const}}
\def\ie{{\it i.e.\ \/}}
\def\cf{{\it cf.\ \/}} 
\def\eg{{\it e.g.\ \/}} 
\def\R{{\bf R}}
\def\C{{\bf C}}
\def\H{{\bf H}}
\def\Q{{\bf Q}}
\def\Z{{\bf Z}}

\def\cb{\overline{\C P}{}^2}
\def\vol{{\it vol}}
\def\diam{{\it diam}}
\def\length{{\it length}}
\def\area{{\it area}}
\def\sys{{\it sys}}
\def\mass{{\it mass}}
\def\Image{{\it Image}}
\overfullrule=0pt           

\bigskip\noindent \vskip2cm \centerline{{\large Systolic freedom of
orientable manifolds}%
\footnote*{{\it Ann.\ Scient.\ Ec.\ Norm.\ Sup.,\/} 1998 (to appear)}
 } 

\bigskip\noindent \vskip1cm \centerline{Ivan
Babenko and Mikhail Katz} \vskip2cm 

\bigskip\noindent {\bf Abstract:} In 1972, Marcel Berger defined a
metric invariant that captures the `size' of $k$-dimensional homology
of a Riemannian manifold.  This invariant came to be called the
$k$-dimensional {\it systole.\/} He asked if the systoles can be
constrained by the volume, in the spirit of the 1949 theorem of
C. Loewner.  We construct metrics, inspired by M. Gromov's 1993
example, which give a negative answer for large classes of manifolds,
for the product of systoles in a pair of complementary dimensions.  An
obstruction (restriction on $k$ modulo 4) to constructing further
examples by our methods seems to reside in the free part of real Bott
periodicity.  The construction takes place in a split neighborhood of
a suitable $k$-dimensional submanifold whose connected components
(rationally) generate the $k$-dimensional homology group of the
manifold.  Bounded geometry (combined with the coarea inequality)
implies a lower bound for the $k$-systole, while calibration with
support in this neighborhood provides a lower bound for the systole of
the complementary dimension.  In dimension 4 everything reduces to the
case of $S^2\times S^2$.

\bigskip\noindent {\large 1. Introduction} \medskip\noindent In 1972,
Marcel Berger defined metric invariants that capture the `size' of
$k$-dimensional homology of a manifold.  He asked if these invariants
can be constrained by the volume (see inequality (1.2) and the
historical section 1.1 below).  We construct metrics, inspired by M.
Gromov's 1993 example, which give a negative answer for large classes
of manifolds.  An obstruction to constructing further examples by our
methods seems to reside in the free part of Bott periodicity (\cf
Theorem 3 below and Lemma 5.5). 

We show that every compact orientable manifold $X$ of dimension $n\geq
3$ admits metrics $g$ of arbitrarily small volume satisfying the
inequality
$$vol_{n-1}(M)\ \length(C)>1 \eqno{(1.1)} $$ for every noncontractible
closed curve $C\subset X$ and every orientable non-separating
hypersurface $M\subset X$ (the result is meaningful whenever the first
Betti number of $X$ is nonzero).  Further results of this type appear
below (the theorems are stated in section 1.6).

We discuss the following special cases: projective spaces (Example 5.3
and Remark 5.6), the 3-torus (Example 3.5), products of spheres
(Proposition 4.2).

Note that for $n= 2$ such metrics do not exist.
For example, every metric 2-torus $T^2$
satisfies 
$${\length^2(C)\over \area(T^2)}\leq {2\over\sqrt{3}},
\eqno{(1.2)}
$$
where $C$ is its shortest noncontractible loop,
by Loewner's theorem (\cf [34], p.\ 295-296 and [43], [6], [7]).

The invariant that captures the notion of the `size' of
$k$-dimensional homology of $X$ is the systole of M. Berger [6].  We
define the $k$-systole, $\sys_k(g)$, of a Riemannian manifold $(X,g)$
as the infimum
$$\sys_k(g)=\inf_M\vol_k(M)
\eqno{(1.3)}
$$ over all nonnullhomologous $k$-cycles $M$ represented by maps of
manifolds into $X$ (\cf Lemma 4.1).

We also define the $k$-mass, $\mass_k(\alpha)$, of a $k$-dimensional
homology class $\alpha$ as follows:
$$\mass_k(\alpha)=\lim_{i\rightarrow\infty}{1\over i}\inf_{M\in
i\alpha}\vol_k(M)
\eqno{(1.4)}
$$ where the infimum is taken over submanifolds $M$ representing
$i\alpha$.  The $k$-mass of $(X,g)$ is the infimum of masses of
nontorsion $k$-dimensional integer homology classes (often called the
{\it stable} $k$-systole).  

\medskip\noindent {\bf 1.1 The origin of the problem.}  Following the
pioneering work of C. Loewner and P. Pu [43] (and, later, R. Accola,
C. Blatter), M. Berger [6] defined the systolic invariants in 1972 in
the pages of the Annales Scientifiques de l'Ecole Normale
Sup\'erieure, and asked if they can be constrained by the volume (see
the excellent survey [8], as well as [10], section TOP.1.E, page 123).

M. Gromov [31] proved in 1983 that this is the case for the top power
of the homotopy 1-systole, for a large class of manifolds called {\it
essential} (see below), which includes aspherical manifolds as well as
real projective spaces.  Namely, for these manifolds one has
$$(\pi \sys_1(g))^n\leq  \Const_n\ \vol(g)
\eqno{(1.5)}
$$ for a positive $ \Const_n>0$ depending only on $n=\dim(X)$, where
$\pi \sys_1(g)$ is the length of the shortest noncontractible loop for
the metric $g$ on $X$.  

An oriented manifold $X$ is {\it inessential} if the inclusion of $X$
in the classifying space $B\pi_1=K(\pi_1(X),1)$ admits a retraction
(fixing the 1-skeleton) to the $(n-1)$-skeleton of $B\pi_1$.  Then the
1-systole is not constrained by the volume, as shown by the first
author in [1], p.\ 34.  We will apply a similar technique which
involves maps to complexes obtained from $X$ by attaching cells of
dimension at most $n-1$ (\cf section 6).

\medskip\noindent {\bf 1.2 Is volume a constraint?}  So far we have
mentioned no positive results on constraining the higher dimensional
systoles $\sys_k$ ($k\geq 2$) by the volume.  Such inequalities do
exist for the $k$-mass.  Indeed, the product of {\it masses} in
complementary dimensions {\it is} constrained by the volume:
$$\mass_k(g)\ \mass_{n-k}(g)\leq{\rm Constant}\ vol(g)
\eqno{(1.6)}
$$ for a suitable Constant depending only on the topology of $X$, for
a wide class of manifolds (\cf [30], p.\ 60; or [35], sections
4.36-38; and [36]).

As late as 1992, Gromov was optimistic about the existence of similar
inequalities for systoles.  Thus at the end of section $4.A_3$ of [32]
one finds a discussion of the systolic $(k,n-k)$ inequality, and a
program of study in terms of conditions on geometry and curvature.
Section $4.A_5$ of [32] contains the following statement: ``In
general, we conjecture that all non-trivial intersystolic inequalities
for simply connected manifolds are associated to multiplicative
relations in the cohomology in the corresponding dimensions. This
conjecture applies to the mass as well as to the volume and moreover
it should be refined in the case of mass [...]''  (the text goes on to
propose such a refinement). The version of [32] published in 1996
modifies this sentence (see [34], p.\ 354) in the light of Gromov's
subsequent developments.  

\medskip\noindent {\bf 1.3 The first counterexample.}  It was invented
in 1993 by M. Gromov on $X=S^3\times S^1$ and described by M. Berger
in [8], p.\ 301 (see also [35], section 4.45 and [38], section 1).
The same paper [8] states Gromov's results concerning the existence of
further examples on products of spheres.  The second author wrote to
M. Berger in February, 1994, asking about the details of Gromov's
construction.  M. Berger immediately replied [9] that he was
awaiting a ``new manuscript from Gromov, or now from you!''  It is the
warm encouragement of M. Berger that initially stimulated the second
author's interest in the problem.

With the emergence of Gromov's 1993 example, the situation changed
drastically.  Most inequalities involving the higher systoles are now
conjectured to fail.  

The inequalities tend to be violated even among homogeneous metrics,
as the 1993 example shows.  It would be interesting to explore the
Lie-theoretic origin of the $S^3\times S^1$ example and its
generalizations.  A homogeneous example in middle dimension on a
product of spheres appears in [5] and [35], Appendix D, Theorem D.3.

\medskip\noindent {\bf 1.4 `Massive' metrics and instability.}  To
present a different point of view on Gromov's construction, it is
convenient to introduce the following terminology.  A sequence $g_j$
of metrics is called {\it massive\/} if each $g_j$ satisfies equality
in the mass inequality $(1.6)$ (for a fixed Constant
independent of $j$); it is called $k$-{\it unstable\/} if
$$\lim_{j\rightarrow\infty}{sys_k(g_j)\over mass_k(g_j)}=\infty\ .
\eqno{(1.7)}
$$  In this language, the 1993 metric is massive
and 1-unstable (\cf Remark 3.2).

Note that $\mass_{n-1}=sys_{n-1}$ by [23], statement 5.10, p.394 
(\cf [30], p.\ 59).

Its 1-instability is due to the fact that at the level of the
universal cover, the step of the covering transformation $T$ is far
greater than the displacement of the hypersurface $S^3$ (\ie the
distance between $S^3$ and $T(S^3)$), due to rapid Hopf rotation.
Gromov privately speculated that the same could be achieved for
$X=S^2\times S^1$ despite the absence of a free circle action, by
stretching in the neighborhood of the fixed points.  Such stretching
was formalized by L. Berard Bergery and the second author in [4] in
terms of the nilmanifold $N$ of the Heisenberg group.

\medskip\noindent {\bf 1.5 Semidirect products and Sol Geometry.}
Note also that the fundamental group of $N$ is a semidirect product.
C.~Pittet [42] clarified the role of semidirect products in the
construction of the example on $S^2\times S^1$, reformulating it in
terms of Sol geometry (see Remark 2.2).  He generalized this example
to manifolds of the form $M\times S^1$ satisfying $H_1(M)=0$ by
localizing the construction of [4]. 

\medskip\noindent {\bf 1.6 Statement of results.}  There are two
rather different systolic problems: in middle dimension, on the one
hand, and in a pair of complementary dimensions, on the other.  The
case of middle-dimensional freedom is treated in [2].  In the
present paper we are concerned with the case of a pair of distinct
complementary dimensions, and also the special case $n=4$ (see Theorem
4 below and Remark 3.6).

\medskip\noindent
{\bf Definition.} We say that $X$ is systolically
{\it `$(k,n-k)$-free'\/} if
$$\inf_g {vol(g)\over \sys_{n-k}(g) \sys_k(g)}=0,
\eqno{(1.8)}$$ 
where the infimum is taken over all Riemannian metrics on $X$. 
If $n=2k$, we say that $X$ is $k$-free.

\medskip\noindent We now state our main results.  Most of the results
in this area are due to M. Gromov.  As described above, the first
results on freedom published with detailed proofs are by L. Berard
Bergery and M. Katz [4], M. Katz [38], and C. Pittet [42].

\medskip\noindent 
{\bf Theorem 1.}
Let $X$ be a compact orientable $n$-manifold, $n\geq 3$.  Then $X$ is
systolically $(1,n-1)$-free.  

\medskip\noindent {\bf Theorem 2.}  Let $X$ be a compact $n$-manifold,
$n\geq 5$, whose fundamental group is free abelian.  Then $X$ is
systolically $(2,n-2)$-free.  

\medskip\noindent {\bf Theorem 3.} Let $X$ be a compact orientable
$(k-1)$-connected $n$-manifold, where $k<{n\over 2}$.  Then $X$ is
systolically $(k,n-k)$-free if $k$ is not divisible by 4.  If $k=4$,
then $X$ is $(4,n-4)$ free provided that $p_1(X)=0$.

\medskip\noindent The 2-freedom of a simply
connected 4-manifold $X$ can be characterized as follows: $X$ admits
metrics of arbitrarily small volume such that every noncontractible
surface inside it has at least unit area.  \medskip\noindent {\bf
Theorem 4.} The following three assertions are equivalent:

(i) the complex projective plane $\C P^2$ is 2-free;

(ii) $S^2\times S^2$ is 2-free;

(iii) every simply connected 4-manifold is 2-free.  \medskip\noindent
{\bf 1.7 The construction.}  To explain the idea in more detail, let $
A\subset X$ be a $k$-dimensional submanifold whose connected
components (rationally) generate $H_k(X)$.  Assume that the normal
bundle of $A$ is trivial, so that its tubular neighborhood in $X$ is
diffeomorphic to $A\times B^{n-k}$ (\cf [14]).  Let $R\subset B^{n-k}$
be a codimension 2 submanifold with trivial normal bundle (for
example, an $(n-k-2)$-sphere).  Then the boundary of a tubular
neighborhood of $R\subset B^{n-k}$ is diffeomorphic to $R\times T^1$.
Assume furthermore that $A$ splits off a circle, \ie $A=C\times S$
where $C$ is the circle.  Let $\Sigma\subset X$ be the hypersurface
$$\Sigma=A\times R\times T^1=T^2\times L$$ where $T^2=C\times
T^1$ is the 2-torus and $L=S\times R$.  A tubular neighborhood of
$\Sigma\subset X$ is a cylinder $\Sigma\times I=Y\times L$ where
$Y=T^2\times I$ is a cylinder on the 2-torus.  We construct direct sum
metrics on $\Sigma=Y\times L$ which are fixed on $L$.  Meanwhile, we
will describe a sequence of special metrics on $Y$ in Lemma
2.1, which give rise to the free metrics on $X$.

Theorem 3 may be viewed as the most general result one can obtain by
generalizing the construction of free metrics on products of spheres
(see Proposition 4.2).  One may ask if its modulo 4 hypothesis can be
removed.  While this may turn out to be possible, the free metrics on
products of spheres would probably still be the starting point.  This
paper may therefore be viewed as an effort to understand the domain of
applicability of the construction described in the previous paragraph,
starting with a codimension 2 submanifold in the fiber of the
(trivial) normal bundle of a split submanifold in X whose connected
components generate $H_k(X)$.

Similar results for the spectrum of the Laplacian have appeared in
[20], [17].

\medskip\noindent {\bf 1.8 Expansion-contraction and Thom's
theorem.}  From the point of view of Thom's theorem, the construction
can be described as follows.  Following Gromov ([34], section
$4.A_5$), we use Thom's theorem to choose a submanifold below the
middle dimension satisfying the following two properties: (a) its
normal bundle is trivial; (b) its connected components generate
(rationally) the homology group whose classes they represent.  The
construction takes place within a fixed trivialized tubular
neighborhood of the submanifold.  Gromov described an
expansion-contraction procedure (\cf [33], section 2) with the
desired effect on the systole of the complementary dimension.  We
combine the expansion-contraction with a volume-saving twist.

In more detail, Thom's theorem allows us to carry out the
expansion-contraction procedure in a nieghborhood of a suitable
$k$-dimensional submanifold described above.  The resulting metrics
have large $(n-k)$-dimensional systole (compared to the volume).  Here
the $(k,n-k)$ systolic inequality is not violated, as the volume is
too large.  To decrease the volume, we introduce a nondiagonal
coefficient in the matrix of the Riemannian metric.  This requires a
splitting of both the $k$-dimensional class and its tubular
neighborhood.  It is at this point that the modulo 4 condition on $k$
comes in.

The construction works whenever a suitable multiple of each
$k$-dimensional homology class contains a representative $A$ with the
following two properties: (a) the normal bundle of $A$ is trivial; (b)
$A$ splits off a circle factor in a Cartesian product.  It suffices to
find a representative which is a sphere with trivial normal bundle.
This is done using Thom's theorem and Bott periodicity.

We construct a sequence of metrics $g_j$ for which the quotient (1.8)
becomes smaller and smaller.  The metric $g_j$ contains two regions:
one where the geometry is fixed (\ie independent of $j$), and another
(locally isometric to the left-invariant metric on the Heisenberg group)
where the geometry is `periodic' with $j$ periods.  The bounded
geometry implies a uniform lower bound for the $k$-systole.  Our key
technical tools here are the coarea inequality and the isoperimetric
inequality of Federer and Fleming.  A suitable calibrating form with
support in the `periodic' region shows that the $(n-k)$-systole grows
faster than the volume.

\medskip\noindent In section 2, we describe a local version of the
construction of free metrics used in [4].  We exploit it in section 3
to prove Theorem 1.  We generalize the construction to higher
$k$-systoles in section 4 for products of spheres and in section 5 in
the general case.  Theorem 4 is proved in section 6, using Whitehead
products and pullback arguments (\cf [1]) for metrics.

\medskip\noindent {\bf 1.9 Acknowledgment.} The second author is
grateful to M. Berger for his warm encouragement and interest in the
progress of the systolic problem.  Thanks are also due to T.~Hangan
for providing a helpful calculation on the Heisenberg group.

\bigskip\noindent {\large 2. Calibration, bounded geometry, and the
Heisenberg group}

\medskip\noindent We present a local version of the
construction of [4] in Lemma 2.1 below.  The significance of the lemma
is that the product of the 2-torus and the interval admits metrics
with enough room for cylinders of unexpectedly large area (but see
(2.9)).

To distinguish two circles which play different roles in the
construction, we denote them, respectively, $T^1$ and $C$.  Consider
the cylinder $M=T^1\times I$, circle $C$, and the manifold
$$Y=C\times M=T^2\times I.$$ We view $M$ as a relative cycle in
$H_2(Y,\partial Y)$.  The 2-mass of the class $[M]$ is the infimum of
areas of {\it rational} cycles representing it, where
$\area(\Sigma_i(r_i \sigma_i))=\Sigma_i r_i \area(\sigma_i)$, $r_i\in
\Q$.

\medskip\noindent {\bf Lemma 2.1.}  
The manifold $Y=T^2\times I$ satisfies
$$\inf_g{\vol(g) \over \sys_1(g)\mass_2([M])} = 0,$$ where the infimum
is taken over all metrics whose restriction to each component of the
boundary $\partial Y=T^2 \times \partial I$ is the standard `unit
square' torus satisfying $\length(C)=\length(T^1)=1$.  

\medskip\noindent {\it Proof.}  More precisely, we
will show that there exists a sequence of metrics $Y_j=(Y,g_j)$
satisfying the following four conditions:

(i) the restriction of $g_j$ to the boundary $T^2\times\partial I$ at
each endpoint is the standard unit square metric for which $T^1$ and
$C$ have unit length;

(ii) the 1-systole of $g_j$  is uniformly bounded from below;

(iii) the volume of $g_j$ grows at most linearly in $j$;

(iv) the 2-mass of $[M]\in H_2(Y,\partial Y)$ grows at least
quadratically in $j$.

We present a shortcut to the explicit formula for the solution.  A
reader interested in the method of arriving at such a formula
geometrically can consult [4], p.\ 630.  Consider the following metric
$h(x)$ in the $y,z$-plane depending on a parameter $x\in \R$:
$$h(x)(y,z)=(xdy-dz)^2+dy^2
\eqno{(2.0)}
$$ Let $I=[0,2j]$.  For $x\in I$, set
$\widehat{x}=j-|x-j|=\min(x,2j-x)$.  Consider the fundamental domain
$D= \{0\leq x\leq 2j,\ 0\leq y\leq 1,\ 0\leq z\leq 1\}$. The metric
$$g_j=h(\widehat{x})(y,z)+dx^2,
{\rm\ where\ } (x,y,z)\in D
\eqno{(2.1)}
$$ gives rise to a metric $g_j$ on $T^2\times I$ once we identify
the opposite sides of the unit square in the $yz$-plane.  The circles
$T^1,C \subset T^2=\R^2/\Z^2$ are parametrized, respectively, by the
$y$- and $z$-axis.  The circles $T^1$ and $C$ have lengths,
respectively, $\sqrt{\widehat{x}^2+1}$ and 1 with respect to the
metric $h(\widehat{x})=(1+\widehat{x}^2)dy^2+dz^2-2\widehat{x} dy\,dz$
on $T^2\times \{x\}$.  The length of $T^1$ thus stretches from 1 to
$\sqrt{j^2+1}$ and then shrinks back to 1, so as to satisfy (i).  Our
choice of the metric such that the length of $T^1\subset T^2\times
\{x\}$ is $\sqrt{1+\widehat{x}^2}$ rather than simply $\widehat{x}$
results in local homogeneity (see below).

The metric $g_j$ is symmetric with respect to the midpoint $x=j$. 
The map 
$$\psi: (x,y,z)\mapsto (x+1,y+z,z)\eqno{(2.2)}$$ defines an isometry
from $T^2\times [i-1,i]$ to $T^2\times [i,i+1]$ for $i=1,\ldots,j-1$.
Thus the metric $g_j$ on $T^2\times [0,j]$ is 1-periodic (\cf (2.6)),
proving (ii).

The $x$-axis projection $Y\rightarrow I$ 
is a Riemannian submersion over an interval of length $2j$
with fibers of unit area, hence the volume estimate (iii). 

We obtain the lower bound (iv) for the 2-mass of $[M]$ as in [4], pp.\
625-626 by a calibration argument.  First we calculate the area of the
cylinder $M$.  Recall that the universal cover of $M$ is a subset of
the $xy$-plane in our coordinates.  Since $T^1\times \{x\}$ has length
$\sqrt{\widehat{x}^2+1}$, we have $$\area(M)=\int_{0}^{2j}\length(T^1\times
\{x\})dx= 2\int_0^{j}\sqrt{x^2+1}\ dx \sim j^2.\eqno{(2.3)}$$ The
calibrating form is the pullback of the area form of $M$ by the
nearest-point projection.  The projection, while not
distance-decreasing, is area-decreasing.  Consider the 2-form
$$\alpha=\sqrt{1+x^2}dx\wedge d\lambda {\rm\ where\ }
\lambda=y-{x\over 1+x^2}z.\eqno{(2.4)}$$ Note that
$\sqrt{1+x^2}\alpha=*(dz)$ where $*$ is the Hodge star of the
left-invariant metric.  This formula is easily verified with respect
to the orthonormal basis of left-invariant forms $dx$, $dy$, $dz-xdy$
(we owe this remark to T.~Hangan).  (Note that in Gromov's example,
the calibrating form is in fact the Hodge star of the projection to
the circle fibre, suitably normalized).  The restriction of $\alpha$
to $M$ coincides with the area form of $M$ for $x\leq j$, and $\alpha$
is of unit norm (see [4], pp.\ 627-628 for details).  Let $\phi_j(x)$
be a partition of unity type function with support in $]0,j[$ and such
that $\phi_j(x)=1$ for $x\in[1,j-1]$.  The form $\phi_j\alpha$ is
closed.  Let $M'$ be any rational cycle representing the class
$\epsilon[M]\in H_2(Y,\partial Y)$, where $\epsilon=\pm 1$.  As in
[4], p.\ 626, we have from (2.3):
$$\mass(\epsilon[M])=\inf_{M'}\area(M')\geq\int_{M'}\epsilon\phi_j\alpha=
\int_{M}\phi_j\alpha\geq\int_1^{j-1}xdx\sim{j^2},\eqno{(2.5)}$$
proving Lemma 2.1.  \medskip\noindent {\bf Remark 2.2.} The
isometry $\psi:(x,y,z)\mapsto (x+1,y+z,z)$ acts in the $yz$-plane by
the unipotent matrix $\left[\matrix{1&1\cr 0 &1\cr}\right]$.  Pittet
[42] replaced this action by a hyperbolic one, producing relative
2-cycles with exponential growth of the 2-mass, rather than quadratic
as in (iv) above.  \medskip\noindent {\bf Remark 2.3.} The
universal cover of the `half' $T^2\times [0,j]$ is isometric to the
subset defined by the condition $0\leq x\leq j$ of the Heisenberg
group of unipotent matrices $$\left[ \matrix{1&x&z\cr 0&1&y\cr
0&0&1\cr}\right];\ x,y,z\in\R$$ with the standard left-invariant
metric $dx^2+dy^2+(dz-xdy)^2$ (\cf [18], p.\ 67; [26], p.\ 227).
The isometry $\psi$ is left multiplication by the matrix
$$\left[\matrix{1&1&0\cr 0&1&0\cr 0&0&1\cr}\right].$$ Let $N=G/\Gamma$
be the standard nilmanifold of the Heisenberg group $G$, where
$\Gamma$ consists of matrices with integer entries.  Factoring by the
iterates of $\psi$, we obtain a projection $f: T^2\times
[0,j]\rightarrow N$ which is a local isometry.  Now take the interval
$I=[0,2j]$ and fold it in two at $x=j$, \ie send $x$ to
$\min(x,2j-x)$.  Let $g: Y_j\rightarrow T^2\times [0,j]$ be the
resulting folding map on $Y_j$.  The {\it distance decreasing}
projection $f\circ g:Y_j\rightarrow N$ induces a monomorphism at the
level of the fundamental groups.  Therefore $$\sys_1(Y_j)\geq \pi
\sys_1(N), \eqno{(2.6)}$$ where $\pi \sys_1$ is the length of the
shortest noncontractible loop.  Meanwhile, the induced homomorphism in
1-dimensional homology sends the class $[C]$ to 0, as the $z$-axis in
$G$ projects to a loop in the center of the fundamental group $\Gamma$
of $N$, which is also its commutator subgroup.  

\medskip\noindent {\bf Remark 2.4.} As an abstract group, $\Gamma$ is
represented by generators $x,y,z$ and relations $[x,y]=z$, $[x,z]=1$,
$[y,z]=1$, where $[a,b] =aba^{-1}b^{-1}$.  Let $a^{(b)}=bab^{-1}$.
Note that for every positive integer $j$ we have the following
relation in $\Gamma$:
$$z^j  = y^{(x^j)}y^{-1}, \eqno{(2.7)}$$ which is the combinatorial
antecedent of (2.8) below.  Indeed, $z^j y =
z^{j-1}xyx^{-1}y^{-1}y=z^{j-1}y^{(x)}=(z^{j-1}y)^{(x)}$ since $z$ is
central.  In the same way, $z^j y =
(z^{j-2}y)^{(x^2)}=\ldots=y^{(x^j)}$.  

Let $$v=\left[\matrix{1\cr j\cr}\right],\ A=\left[\matrix{1+j^2&-j\cr
-j &1\cr}\right],$$ so that ${}^tv A v =1$.  It follows that the
shortest loop in the class $[T^1 + jC]\in H_1(T^2\times\{j\})$ has
unit length, since the torus $T^2\times\{j\}$ is equipped with the
metric $dy^2+(dz-j\,dy)^2$.  Thus by sliding the curve $C\subset T^2$
to the value $x=j$, we obtain $$\mass_1[T^1+jC]\leq 1.\eqno{(2.8)}$$
\medskip\noindent {\bf Remark 2.5.} The metric $h_x$, defined by the
matrix $\left[\matrix{1+\widehat{x}^2&-\widehat{x}\cr -\widehat{x}
&1\cr}\right]$, where $\widehat{x}=\min(x,2j-x)$, is flat of unit
area.  The flat tori of unit area and unit 1-systole used in our
construction all lie in a compact part of the moduli space of tori,
namely the interval $s+\sqrt{-1}$ for $s\in [-{1\over 2},{1\over 2}]$
in the standard fundamental domain in $\C$.  The diameter of each of
these tori is less than 1, so that our manifold $Y_j$ is rather
narrow: $$diam_1(Y_j)<1,\eqno{(2.9)}$$ where the 1-diameter of a
Riemannian manifold $X$, denoted $\diam\sb 1\,X$, is the infimum of
numbers $\varepsilon>0$ such that there exists a continuous map from
$X$ to a graph with the property that the inverse image of every point
has diameter $\le\varepsilon$.

More precisely, given a map $f:X\rightarrow \gamma$, from a Riemannian
manifold $X$ to a graph $\gamma$, define the size $s(f)$ by
$s(f)=\sup_{x\in\gamma}\diam(f^{-1}(x)).$ Then the 1-diameter of $X$
is the least size of a map from $X$ to a graph: $\diam_1(X)=
\inf_{\gamma,f}s(f)$, where the infimum is taken over all graphs
$\gamma$ and all continuous maps $f:X\rightarrow\gamma$.  Setting
$\gamma=[0,2j]$ and $f=$ the $x$-coordinate, we obtain (2.9).
\medskip\noindent {\bf Remark 2.6.} The construction of the manifold
$Y$ of Lemma 2.1 can be summarized as follows.  One takes a torus
bundle, say $N$, over a circle, whose glueing automorphism $A$ is not
an isometry.  Let $v$ be a lattice vector whose images increase
indefinitely in length under the iterates of $A^{-1}$.  One chooses a
fixed metric on $N$ and pulls it back to the infinite 3-dimensional
cylinder $Y_\infty=T^2\times \R$, where $\R$ is the universal cover of
the circle.  This produces a periodic metric on $Y_\infty$.  Consider
the 2-dimensional cylinder $M\subset Y_\infty$ which is the circle
subbundle spanned by the direction $v$.  The key point now is that the
area of $M$ grows faster than the volume of $Y_\infty$ due to the
choice of the direction $v$.  One then takes a long piece of
$Y_\infty$ and doubles it as described above to obtain the desired
metric.

\medskip\noindent {\bf Remark 2.7.} A.  Besicowitch
[11] in 1952 exhibited a different type of (1,2)-freedom on a
3-dimensional manifold with boundary, namely the cylinder $D^2\times
[0,1]$ (\cf [13], p.\ 296), disproving a conjecture of Loewner's.

\bigskip\noindent {\large 3. Construction in dimension/codimension
one} \medskip\noindent {\bf Proposition 3.1.} Every oriented
$n$-manifold with $b_1(X)=1$ is systolically $(1,n-1)$-free if $n\geq
3$.  \medskip\noindent {\it Proof.} The construction of free metrics
is local in a neighborhood of a loop $C$ which generates $H_1(X,\Z)$
modulo torsion.  Since $X$ is orientable, the normal bundle of $C$ is
trivial, and its tubular neighborhood in $X$ is diffeomorphic to $
C\times B^{n-1}$.  Let $L\subset B^{n-1}$ be a codimension 2
submanifold which then has trivial normal bundle (for example, an
$(n-3)$-sphere).  The boundary of a tubular neighborhood of $L\subset
B^{n-1}$ is diffeomorphic to $L\times T^1$.  Let $\Sigma\subset X$ be
the hypersurface $$\Sigma= C\times L\times T^1=T^2\times L$$ where
$T^2=C\times T^1$ is the 2-torus.  A tubular neighborhood of
$\Sigma\subset X$ is a cylinder $\Sigma\times I=Y\times L$ where
$Y=T^2\times I$ is a cylinder on the 2-torus.  We construct direct sum
metrics on $Y\times L$ which are fixed on $L$.  The special metrics on
$Y$ of Lemma 2.1 give rise to the free metrics on $X$.

A similar technique was used by C. Pittet [42].  His idea was to use
the torus $\Sigma=T^{n-1}$ (where we use $T^2\times L$), but only 2
circles are actually needed.  

What happens metrically can be described
as follows.  We choose a metric on $X$ which is a direct sum in

$$\Sigma\times I=T^2\times L \times I,\eqno{(3.1)} $$ where $T^2$ is
the standard unit square torus (with $C$ and $T^1$ both of unit
length) and $L$ has unit volume.  We now modify the metric in
$\Sigma\times I=Y\times L$ by means of the metric $Y=(Y,g_j)$
of Lemma 2.1, while $L$ keeps the same fixed metric of unit
volume, and the metric on $Y_j\times L$ is a direct sum.  Condition
(i) of Lemma 2.1 ensures that the metric varies continuously across
$\partial I$.  Denote the resulting Riemannian manifold $X_j$.

The metric of $X_j$ stays the same on the complement of $\Sigma\times
I$, while on the region $Y_j\times L$ it is periodic in the sense
given in the proof of Lemma 2.1.  A loop of length less than 1 is
contained either in the cylinder, or in the 1-neighborhood of its
complement.  In the latter case, it can be viewed as a loop in $X_j$
for $j=1$.  Thus $$\sys_1(X_j)\geq\min(\sys_1(X_1),\pi
\sys_1(N))\eqno{(3.4)}$$ by (2.6).  Hence the 1-systole (as well as
the homotopy 1-systole, justifying the first inequality (1.1) of the
paper) is uniformly bounded from below as $j$ increases.

The calibrating $(n-1)$-form $\beta$, supported on $Y\times L$, is
obtained from the form (2.4) of Lemma 2.1 by exterior product with the
volume form $\vol_L$ of $L$:

$$\beta=\phi_j\alpha\wedge \vol_L.  \eqno{(3.7)}$$ 

Choose a hypersurface $\tilde M\subset X$ `dual' to $C$, which has
standard intersection $B^{n-1}$ with a tubular neighborhood of $C$.
Then $\tilde M\cap Y_j =M$, the 2-dimensional cylinder of Lemma 2.1.
Hence

$$\sys_{n-1}(X_j)\geq \mass_{n-1}([\tilde M])\geq
\int_{\tilde M}\beta=vol(L)\int_{M}\phi_j\alpha \sim j^2,
\eqno{(3.9)}$$ proving Proposition 3.1.

\medskip\noindent {\bf
Remark 3.2.} This is a convenient time to explain why these metrics on
X are 1-unstable, \ie the 1-mass tends to 0.  Indeed, by construction
the curve $T^1$ is contractible in $X$.  Hence we have from (2.8),
$$\mass_1[C]={1\over j}\mass_1[jC]={1\over
j}\mass_1[T^1+jC]\leq{1\over j} \rightarrow 0.\eqno{(3.5)}$$ The map
$f\circ g$ of section 2.3 extends to a distance-decreasing map $h$ to
a fixed CW complex $X\cup (N\times L)$, $$h:X_j\rightarrow X\cup
(N\times L), 
\eqno{(3.6)}
$$ where manifolds $X$ and $N\times L$ are glued along the common
hypersurface $T^2\times L$.  The induced homomorphism at the level of
fundamental groups has nontrivial kernel generated by the curve $C$
(the commutator of $T^1$ and the other generator of $\pi_1(N)$).

Note that a contrary, and false, claim was
made in [4], p.\ 626 concerning the homomorphism $p_5$, in order to
prove a uniform lower bound for the 1-systole (see (4.9)).  However,
such a bound follows immediately from the bounded geometry resulting
from the construction, without recourse to the the fixed CW complex.

\medskip\noindent {\bf Remark 3.3.} To a myopic observer, $X_j$
looks like an interval of length $2j$.  Indeed, let $\gamma$ be the
graph defined by the connected components of the level sets of the
distance function from the subset $X_{+}\subset X_j$.  By (2.9),
$$diam_1(X_j)\leq\max(1,diam(X_{+})).\eqno{(3.10)}$$ 

\medskip\noindent {\bf Remark 3.4.}  We explain the choice of
$C$ and $M$ above in terms of Morse theory.
Let $a$ be a
generator of the cohomology group $H^1(X,\Z)=\Z$.  Then $a$ defines a
map $X\rightarrow S^1$ which can be made into a Morse map.  The
inverse image of a regular point is a smooth submanifold $M$.  We cut
$X$ open along $M$ (\ie cut $S^1$ at the regular point) to obtain a
cobordism with top and bottom given by $M$, while the map to $S^1$
becomes an ordinary Morse function on this cobordism.  Next, take the
`same' point at top and bottom (\ie the same as a point of $M$) and
join them in the cobordism by a path with the following two
properties:

(i) it avoids the critical points of the Morse function;

(ii) the Morse function is strictly monotone along the path.

\noindent
Now if we glue top and bottom to form $X$ again,
this path becomes a loop $C$ which meets $M$ transversely
in exactly one point $p\in X$.

\medskip\noindent $\underline{\hbox{{\it Proof of Theorem 1.}}}$ In
the general case $b=b_1\geq 1$, we will first construct curves
$C_1,\ldots,C_b$ and hypersurfaces $M_1,\ldots,M_b$ such that the
intersection $M_i\cap C_k$ contains exactly $\delta_{ik}$ points, as
follows.  We choose a $\Z$-basis $a_1,\ldots,a_{b}$ for the
1-dimensional cohomology $H^1(X,\Z)$.  We define $M_i$ to be the
inverse image of a regular point of a map $X\rightarrow S^1$ defined
by $a_i$.  The $b$ maps define the period map $X\rightarrow T^{b}$
which is surjective at the level of the fundamental groups:
$\pi_1(X)\rightarrow \pi_1(T^{b})=\Z^{b}$ (Poincar\'e duality).
We now choose curves
$C_1,\ldots,C_{b}\subset X$ whose images represent the standard
generators of $\pi_1(T^{b})=\pi_1(S^1)+\ldots+\pi_1(S^1)$.  Then
$a_i(C_k)=\delta_{ik}$ and so the algebraic number of points in the
intersection $M_i\cap C_k$ is exactly $\delta_{ik}$.  

To eliminate points of intersection with negative intersection index,
we choose two {\it adjacent} points on $C_k$ with opposite
intersection indices.  We now perform a surgery on $M_i$ by removing a
little disk around each of the two points and attaching a thin tube to
$M_i$ along the piece of $C_k$ joining the two points.  In this way we
remove all negative intersections.  In particular, we may assume that
the intersection $M_i\cap C_k= \emptyset$ is empty if $i\not=k$.  We
choose a small $\epsilon>0$ so that $M_i\cap(\cup_\epsilon
C_k)=\emptyset$ for all $i\not= k$.  We then insert $b$ copies of
$Y_j\times L$ inside $\cup_\epsilon C_k$ as in (3.1) for
$k=1,\ldots,b$, to obtain a new Riemannian manifold $X_j$
diffeomorphic to $X$.

The calibration argument is generalized as follows.  Let 
$\beta_k=\phi_j\alpha\wedge \vol_L$ be the closed
$(n-1)$-form supported in $\cup_\epsilon C_k$.
Then 
$$\int_{M_i}\beta_k=0 {\rm\ if\ } i\not=k.
\eqno{(3.11)}
$$
Take any nonzero integer class 
$m=\sum_i\epsilon_i d_i[M_i]\in H_{n-1}(X)$
where $\epsilon_i=\pm 1$ and $d_i\geq 0$, $i=1,\ldots,b$.
We use the signs $\epsilon_i$ to specify a calibration form 
$\beta=\sum_k \epsilon_k \beta_k$.  Since the supports
of the $\beta_k$ are disjoint, the form $\beta$ has norm 1.
Let $M'\in m$ be any rational cycle.  Then
$$\vol_{n-1}(M')\geq\int_{M'}\beta =\int_{M'}\sum\epsilon_k \beta_k=
\sum_{i,k} \epsilon_i d_i \int_{M_i}\epsilon_k \beta_k.
\eqno{(3.12)}
$$
In view of (3.11) we have
$$\vol_{n-1}(M')\geq\sum_i d_i \int_{M_1}\beta_1 \geq \int_{M_1}\beta_1
\sim j^2.
\eqno{(3.13)}
$$
This completes the proof of Theorem 1.
\bigskip\noindent
{\bf 3.5. Example of the 3-torus.}
Let $X=T^3=\R^3/\Z^3$.  Define three curves $C_1,C_2,C_3\subset T^3$ 
respectively as the projections of the lines
$\{(t,{1\over 3},{1\over 3})\}$,
$\{({1\over 3},t,{2\over 3})\}$, and
$\{({2\over 3},{2\over 3},t)\}$, where $t\in\R$.
Let $M_i\subset T^3$ be the 2-torus which is the
projection of the coordinate plane perpendicular
to $C_i$, \eg $M_1$ is the projection of $\{(0,s,t); s,t\in\R\}$.
Then $M_i\cap C_k$ consists of $\delta_{ik}$ points.
The boundary of the ${1 \over 9}$-neighborhood of $C_i$
is a 2-torus 
$$T^2_i=\partial(\cup_{1\over 9}C_i)= C_i\times T^1.
\eqno{(3.14)}
$$
We modify the flat metric in a 
${1\over 27}$-neighborhood of each $T^2_i$
to make it into the direct sum of circles
$C_i$ and $T^1$ of unit length.
We cut $X$ open along each $T^2_i$ and insert
three copies of $Y_j$ of Lemma 2.1 along the cuts.
The resulting manifold $(T^3,g_j)$ has 1-systole uniformly bounded
below as $j$ increases, 2-systole growing as $j^2$, and volume as $j$.
This metric bears a formal resemblance to the `Hedlund'
metrics described by V. Bangert [3], p. 278, though it is not obtained
by a conformal coordinate change from the flat metric.

\medskip\noindent{\bf Remark 3.6.}  Taking a product of $(T^3,g_j)$
with a circle of length $j^2$, we obtain metrics on the 4-torus with
2-systole growing faster than the square root of the volume.  In other
words, $T^4$ admits metrics of arbitrarily small volume such that 
every surface inside it representing a nonzero class in
$H_2(T^4)$ has at least unit area.

\bigskip\noindent {\large 4. Products of spheres, coarea inequality, and
intersection number}

\medskip\noindent The product of the
$(1,n-1)$-free metrics of section 3 on an $n$-manifold $X$ with the
sphere $S^{k-1}$ yields $(k,n-1)$-free metrics on $X\times S^{k-1}$.
As we will see, it turns out that it is sufficient to have a product
structure at the level of a $k$-dimensional class containing a
representative $A$ with trivial normal bundle.  If such an $A$ splits
off a circle $C$:
$$A = B \times C\eqno{(4.1)}$$ for some $(k-1)$-dimensional $B$, we
can start pasting in the special metrics of Lemma 2.1.  The idea of
the construction is to keep the factor $B$ of $A$ as a direct summand,
while inserting, as in the case $k=1$, a cylinder obtained by doubling
a piece of the Heisenberg group of length $j\rightarrow\infty$.  Here
the circle $C$ plays the same role as in the local construction of
Lemma 2.1 on $Y_j = T^2 \times I$.  We will illustrate the
construction for the product of two spheres (Proposition 4.2), and
treat the general case in the next section.

All of our lower bounds for the $k$-systole for $k \geq 2$ rely upon a
technique which combines the coarea inequality and the existence of an
intersection number of cycles (dual to the cup product of the dual
classes).  We now present the relevant lemma.  All manifolds are
assumed orientable.

\medskip\noindent {\bf Lemma 4.1.}  Let $D$, $E$, and $G$ be
submanifolds of a manifold $X$ meeting transversely in a single point,
and consider the transverse intersection
$$F = D\cap E.
\eqno{(4.2)}
$$ Let $d\in [D]$ be a cycle defined by the map of a
manifold into $X$.  Then for a dense open set of such maps, the
intersection $f= d\cap E$ is a cycle but not a boundary.

\medskip\noindent{\it Proof.\/} If $D$ and $E$ have complementary
dimensions, the lemma is immediate from the existence of the
intersection number of two cycles, Lemma 10 of J. Schwarz [44] p.\
31 (\cf [27], [39], [22], vol.\ 2, 201O, p.\ 771).  In the general
case, the fact that the class of the intersection of the
representative cycles is independent of the representatives seems to
be difficult to find in the literature.  Note that by the Jiggling
Lemma ([44]. p.\ 24), the intersection $d\cap E$ may be assumed to
be the image of a manifold, as well.  In particular, $d\cap E$ is a
cycle.  By the associativity of set-theoretic intersection, we have
$$(d\cap E)\cap G = d\cap(E\cap G).
\eqno{(4.3)}
$$ Applying the above lemma on the
intersection number twice, we obtain
$$[f] \cdot [G] = [f\cap G] = [d\cap(E\cap G)] = [d] \cdot [E\cap G]
\eqno{(4.4)}
$$
and therefore
$$[f] \cdot [G] = [D] \cdot [E\cap G] = [D\cap E\cap G] = 1\in H_0(X).
\eqno{(4.5)}
$$ Thus $[f]$ is nonzero.

Note that we need to work with maps from manifolds rather than
arbitrary cycles to be able to apply the Jiggling lemma.  Note also
that we do not prove that $[f] = [F]$, only that the class $[f]$ is
nonzero.  Of course, if there is no torsion we can get the stronger
result.

\medskip\noindent {\bf Proposition 4.2.} The manifold $X=S^k\times
S^{n-k}$ is $(k,n-k)$-free for all $n\geq 3$ and all $k$, except
possibly $S^2\times S^2$.

\medskip\noindent {\it Proof.}  For $k=1$, the proposition is a
special case of Theorem 1.  The $k$-freedom of $S^k\times S^k$ for $k
\geq 3$ was established in [38].  We may thus assume that $2\leq
k<{n\over 2}$.  Let $A\subset X$ be a copy of the sphere $S^k$.  We
represent the class $[A]\in H_k(X)$ by an imbedded product
$S^{k-1}\times C\subset X$, where $C$ is a circle (\cf [28], p.\ 33).

This can be done inside an imbedding of $A\times I$ in $X$, where $I$
is an interval.  Join two points of $A\subset A\times I$ by a path
disjoint from $A$.  Remove little $\epsilon$-disks around its
endpoints in $A$.  Attach to $A$ the boundary of the tubular
$\epsilon$-neighborhood of the path.  The resulting hypersurface
$S^{k-1}\times C\subset A\times I$ has trivial normal bundle in
$A\times I$ and hence also in $X$.  Hence its tubular
$\epsilon'$-neighborhood $\cup_{\epsilon'}(S^{k-1}\times C)$ is
diffeomorphic to $S^{k-1}\times C\times B^{n-k}$.  Let $M\subset X$ be
a copy of the sphere $S^{n-k}$ which meets $S^{k-1}\times C$ in a
single point.  We can assume that the intersections are standard:
$M\cap(\cup_{\epsilon'}(S^{k-1}\times C))=B^{n-k}$.  Let $T^1\subset
B^{n-k}$ be an imbedded circle.  For sufficiently small $\epsilon''$,
the boundary of the tubular $\epsilon''$-neighborhood of $T^1$ is
diffeomorphic to $T^1\times S^{n-k-2}\subset B^{n-k}\subset M$.

The hypersurface
$$\Sigma=S^{k-1}\times C\times T^1\times S^{n-k-2}\subset
\cup_{\epsilon'}(S^{k-1}\times C)\eqno{(4.6)}$$ separates $X$ into two
connected components, $X_- = S^{k-1}\times C\times T^1\times
B^{n-k-1}$ and $X_+$ (for which no product structure is available).
Let $T^2=C\times T^1$ and $L=S^{k-1} \times S^{n-k-2}$.
A tubular neighborhood of $\Sigma$ is of the form
$$\Sigma\times I=Y\times L
\eqno{(4.7)}
$$ where $Y=T^2\times I$.  We choose a fixed metric on $X$ outside the
tubular neighborood, in such a way that $T^2$ becomes the unit square
torus (with both $C$ and $T^1$ of unit length).  We then use the
metrics $g_j$ of Lemma 2.1 on $Y$, while the product $Y\times L$ has a
direct sum metric.  We thus obtain the manifold
$$X_j = X_{-} \cup (Y_j\times L)\cup X_{+} \eqno{(4.8)}$$
diffeomorphic to $X$.

The quadratic (in $j$) lower bound for the $(n-k)$-systole of $X_j$ is
obtained by calibration using the form $\beta=\phi_j\alpha\wedge
\vol_{S^{n-k-2}}$ (\cf formula (3.7)).  Proposition 4.2 now follows
from the following lemma.

\medskip\noindent {\bf Lemma 4.3.}  The $k$-systole of $X_j$ is
uniformly bounded from below in $j$.

\medskip\noindent{\it Proof.}  Recall that the metrics $(Y,g_j)$ of
Lemma 2.1 are periodic, while the metric on $X$ outside $Y\times L$ is
fixed.  In the case $k=1$ this was sufficient to obtain a lower bound
for the 1-systole.  However, for $k\geq 2$, the diameter of a small
$k$-cycle is not necessarily small.  We therefore have to give some
proofs using the coarea (Eilenberg's) inequality.  Let $z$ be a
$k$-cycle in $X_j$ representing a nonzero homology class.  

\medskip\noindent {\bf Case 1.}  If $z$ lies in $X_-\cup X_+\subset
X_j$ then it may be viewed as a cycle in $(X,g)$.  Hence $vol(z)\geq
\sys_k(g)$.

\medskip\noindent {\bf Case 2.}  Suppose $z$ avoids $X_{-}\cup X_{+}$,
then $z\subset Y_j \times L\subset X_j$.  Let $\pi: Y_j \times
L\rightarrow S^{k-1}$ be the projection to the first factor of $L$.
Since we have a projection of a Riemannian product, the coarea
inequality applies.  Assuming transversality, we can write
$$vol(z)\geq \int_{S^{k-1}}{\length}(z\cap\pi^{-1}(t))dt
\geq \int_{S^{k-1}} \sys_1(Y_j) dt \geq 
\int_{S^{k-1}} \pi \sys_1(N) dt,\eqno{(4.9)}
$$ again a uniform lower bound (\cf formula (2.6)).  To justify the
last inequality, we apply a relative version of Lemma 4.1 with
$D=S^{k-1}\times C$, $E= Y_j\times S^{n-k-2}$, and $F=C$.  We
conclude that the cycle $z\cap\pi^{-1}(t)\subset Y_j \times
\{t\}\times S^{n-k-2}$ also represents a nonzero multiple of $[C]\in
H_1(\Sigma\times I)$ and therefore participates in the evaluation of
$\sys_1(Y_j \times L)= \sys_1(Y_j)$.

The relevant pairing,
$$H_k(N)\otimes H_{n-k+1}(N,\partial N)\rightarrow
H_1(N)
\eqno{(4.10)}
$$ 
is the homological operation Poincare-Lefschetz dual to the cup
product in cohomology
$$H^{n-k}(N,\partial N)\otimes H^{k-1}(N)\rightarrow
H^{n-1}(N,\partial N),
\eqno{(4.11)}
$$ where $N=\Sigma\times I$.  

\medskip\noindent {\bf General case.} The idea is to cut $z$ at a
narrow place and split it into the sum of 2 cycles which fall into the
2 special cases above.  We argue by contradiction.  Suppose the area
of the $k$-cycle $z_j$ in $X_j$ tends to 0 as $j\rightarrow\infty$.
Let $d$ be the distance function from the subset $X_-\sqcup X_+\subset
X_j$ (on the inserted cylinder, $d$ equals $\min(x, 2j-x)$).  By the
coarea inequality, we find $x_0\in[0,1]$ such that $vol_{k-1}(z\cap
d^{-1}(x_0))\rightarrow 0$.

Note that the
subset $d^{-1}([0,x_0])\subset X_j$ admits a continuous retraction
with a fixed Lipschitz constant to $(X_-\sqcup X_+,g')$, and its
complement, to $\Sigma\times I=Y_j\times L$.

Let $\gamma= z\cap d^{-1}(x_0)$.  Since $vol_{k-1}(\gamma)\rightarrow
0$, it is homologous to 0 in $\Sigma=S^{k-1}\times C\times T^1\times
S^{n-k-2}$ (product of four spheres).  Hence $\gamma$ can be filled in
by a $k$-chain $D$ with $vol_k(D)\rightarrow 0$, so that $\partial D =
\gamma$.  This follows from the isoperimetric inequality for small
cycles in products of spheres, proved in [38], p. 203 (the
generalisation to the case of 4 spheres instead of 2 is
straightforward).  This isoperimetric inequality is an immediate
consequence of the isoperimetric inequality of Federer and Fleming
[24].  It is also a special case of [31], Sublemma $3.4.B'$.  Now we
let $a=(z\cap d^{-1}([0,x_0]))-D$ and $b=z-a$, and apply the
two special cases discussed above, obtaining
$$vol_k(z)\geq\min(vol_k(a),vol_k(b))-vol_k(D) \geq
\min(\sys_k(g),\pi \sys_1(N))-o(1).\eqno{(4.12)} $$ 
This proves Lemma 4.3.

\bigskip\noindent {\large 5. Construction in general codimension and
real Bott periodicity}

\medskip\noindent {\bf Lemma 5.1.} Let $X$ be an $n$-dimensional
orientable manifold and let $k<{n\over 2}$.  Assume that $H_{k-1}(X)$
is torsion-free.  Then $X$ is $(k,n-k)$-free if a suitable multiple of
each $k$-dimensional homology class contains a representative $A$ with
trivial normal bundle such that either $A$ splits off a circle in a
Cartesian product (\ie $A=B\times C$ where $C$ is a circle), or $A$ is
a sphere $S^k$.

\medskip\noindent {\it Proof.}  By Poincar\'e duality
$H_{n-k}(X)=H^k(X)$.  By the universal coefficient formula (\cf [29],
p.\ 194, ex.\ (23.40)), the torsion of this group equals that of
$H_{k-1}(X)$ which vanishes by hypothesis.  Let $M_i$ be an integer
basis for $H_{n-k}(X)$. Let $A_l$, a dual `basis' of $H_k(X)$, in the
sense that the intersection numbers satisfy $M_i \cdot
A_l=\delta_{il}$, \cf section 3.
Recall that the intersection numbers are well defined for each pair 
of cycles of complementary dimensions by [44], p.\ 30-31.

By hypothesis, a suitable multiple of each class $A_l$
contains a representative to which we can apply the construction of
Proposition 3.1.  These multiples are fixed once and for all.
Therefore a systolic lower bound for the multiple implies a lower
bound for the classes $A_l$ themselves.

If $b_k=1$, we carry out the construction of Proposition
3.1, with $L=S^{n-3}$ replaced by $L=B\times S^{n-k-2}$, as in the
proof of Proposition 4.2.  For arbitrary $b_k$ we proceed as in the
proof of Theorem 1.

If $A=S^k$, we add a homologically trivial handle to replace $S^k$ by
$S^{k-1}\times C$, where $C$ is a circle, and argue as before.
\medskip\noindent $\underline{\hbox{{\it Proof of Theorem 2.}}}$ An
oriented rank $n-2$ bundle $\nu$ over an oriented surface $A$ is
determined by the identification of a pair of trivialisations over a
small disk and its complement, along the circle which is their common
boundary, \ie by an element of $\pi_1(SO_{n-2})=\Z_2$ since $n-2\geq
3$.  This $\Z_2$ information is also contained in the second
Stiefel-Whitney class $w_2(\nu)\in H^2(A,\Z_2)=\Z_2$, or equivalently
the second Stiefel-Whitney number $w_2[\nu]=w_2(\nu)[A]\in
H_0(A,\Z_2)= Z_2$ (\cf [41], p.\ 50).

Now let $g:A\rightarrow X$ be an imbedding, and
let $\nu=\nu(A\subset X)$ be its normal bundle in $X$.  Since
$g^*(TX)=TA+\nu$, we have the following identity in $H^2(A)$ by
[41], p.\ 38:
$$g^*w_2(TX)=w_2(TA)+w_1(TA)w_1(\nu)+w_2(\nu)=w_2(\nu)\eqno{(5.1)} $$
since $A$ is spin.  If $X$ is spin, then $\nu$ is automatically
trivial.  Otherwise we assume that $g(A)$ represents an even multiple
of $A_0$, an integer class: $g_*[A]=2[A_0]\in H_2(X)$.  Then
$$w_2[\nu]= w_2(\nu)[A]=g^* w_2(TX)[A]=w_2(TX)(g_*[A])=
w_2(TX)(2[A_0])=0.\eqno{(5.2)} $$ 

\noindent {\bf 5.2 Use of Thom's theorem.}  Every 2-dimensional homology
class in codimension at least 3 can be represented, up to a multiple,
by an imbedded surface $A$ (\cf [46] and [22], p.\ 434).  If
$A\subset X$ has genus $g\geq 2$ then it can be cut along
nullhomologous curves into pieces of genus 1.  Since by hypothesis
$\pi_1(X)$ is abelian, the curves are nullhomotopic.  Hence $A$ is
homologous to a union of tori.  An imbedded 2-sphere in $X$ can also
be turned into a torus by adding a homologically insignificant handle.
Thus every 2-dimensional class can be represented by a union of tori.
Furthermore, the torus admits a self-map of degree 2, which can be
perturbed in this codimension to be made into an imbedding.  The
analysis above now implies that every class in $H_2(X)$ can be
represented, up to a multiple, by a union of imbedded tori with
trivial normal bundles.  We apply Lemma 5.1 to complete the proof of
Theorem 2.

\medskip\noindent {\bf Example 5.3.} The manifold
$X=\C P^{m}$ is $(2,2m-2)$-free for all $m\geq 3$.

\medskip\noindent $\underline{\hbox{{\it Proof of Theorem 3.}}}$ We
first treat two special cases $k=3$ and $k=4$.  

Let $k=3$.  The normal bundle of every imbedded 3-sphere in $X$ is
trivial because $X$ and $S^3$ are oriented so that the normal bundle
is orientable, and an oriented bundle of rank $n-3$ over $S^3$ is
determined by the identification of the trivialisations over two
hemispheres along their common boundary $S^2$, \ie by an element of
$\pi_2(SO_{n-3})$.  The latter group is trivial by a theorem of Elie
Cartan [15], [19].

To give a quick proof of the vanishing of this
group, note that $\pi_2(SO_3)=\pi_2(Spin_3)=0$.  The long exact
sequence of the fibration $SO_m\rightarrow SO_{m+1}\rightarrow S^m$
gives $\pi_2(SO_m)\rightarrow \pi_2(SO_{m+1})\rightarrow \pi_2(S^m)$
and therefore $\pi_2(SO_{n-3})=0$ for all $n\geq 7$.  

Since $X$ is 2-connected, every $3$-dimensional class can be
represented by the map of the 3-sphere, by Hurewicz's theorem (\cf
[25], p.\ 101), and we apply Lemma 5.1 to complete the proof in the
case $k=3$.

\medskip\noindent {\bf Example 5.4.} The manifold $SU(3)$ is
2-connected and therefore (3,5)-free.  Is this still true if the
competing metrics of formula (1.8) are required to be left-invariant?

\medskip\noindent
Let $k=4$ and $p_1(X)=0$.
A calculation similar to (5.1) shows that the first
Pontrjagin class $p_1(\nu)$ of
the normal bundle $\nu$ (of rank $n-4\geq 5$)
of an imbedded 4-sphere vanishes
modulo 2-torsion (\cf [41], Theorem 15.3, p.\ 175).
Since $H^4(S^4)=\Z$ is torsion-free, we have $p_1(\nu)=0$.
We have
$$\pi_4(BSO(n-4))=\pi_3(SO(n-4))=\Z\hbox{{\rm\ for\ }} n\geq 9
\eqno{(5.3)}$$
(\cf [22], vol.\ 4, p.\ 1745, Appendix A, Table 6.VII).
The map from $\pi_3(SO(n-4))=\Z$
to $H^4(S^4)=\Z$, defined by taking the first Pontrjagin class
of the corresponding bundle,
is a nontrivial additive homomorphism
(\cf [41], p.\ 246, 145), hence injective.  It follows
that $\nu$ is the trivial bundle.

Now suppose $k$ is not divisible by 4.  Since $(n-k)-(k-1)\geq 2$, we
have $\pi_{k-1}(SO(n-k))=\pi_{k-1}(SO)$ (the stable group; \cf [21],
chapter 6, paragraph 24, p.\ 233).  By real Bott Periodicity,
$\pi_{k-1}(SO)$ is either 0 or $\Z_2$ if $k$ is not divisible by 4
(\cf [12], p.\ 315).  The theorem now follows from the following
lemma.

\medskip\noindent {\bf Lemma 5.5.}  Let $g:S^k\rightarrow X$ be a
sphere representing an even multiple of a class in $H_k(X)$, where $k$
is not a multiple of 4 and $2k<n$.  Then its normal bundle $\nu$ is
trivial.

\medskip\noindent {\it Proof.}  Recall that for such $k$, real Bott
periodicity implies that the only nontrivial bundle, if it exists, is
distinguished by the $k$-th Stiefel-Whitney class.  Now
$$g^*w_k(TX)=w_k(TS^k)+w_{k-1}(TS^k)w_1(\nu)+\ldots+w_k(\nu)=w_k(\nu)
\eqno{(5.4)}
$$ since the total Stiefel-Whitney class of the sphere is trivial:
$w(S^k)=1$. Writing $g_*[A]=2[A_0]\in H_k(X,Z)$, we obtain 
$$w_k[\nu]= w_k(\nu)[S^k]=g^* w_k(TX)[S^k]=w_k(TX)(g_*[S^k])=
w_k(TX)(2[A_0])=0.  \eqno{(5.5)}
$$
Hence the normal bundle of $g(S^k)$ is trivial.  Lemma 5.5 and Theorem
3 are proved.

\medskip\noindent {\bf Remark 5.6.} The simplest manifold not covered
by Theorem 3 is the quaternionic projective space $\H P^3$.  Its (4,8)
systolic freedom cannot be established by our methods as $p_1(\H
P^3)=4u\in H^4(\H P^3)$ is nonzero, where $u$ is a generator (\cf
[41], p.\ 248).  Note that the systolic 4-freedom of $\H P^2$ follows
by a method similar to Lemma 6.4 below, starting with free metrics on
$S^4\times S^4$, \cf [38] and [2].  We find ourselves in the awkward
situation of being able to establish the systolic freedom of $\C P^3$
(Theorem 2) and $\H P^2$, but not of $\C P^2$ or $\H P^3$.

\bigskip\noindent
{\large 6. Systolic freedom in four dimensions}

\medskip\noindent {\bf Lemma 6.1.}  Let $X$ and $X'$ be manifolds of
the same dimension $n$.  Suppose $X'$ admits a continuous map to the
mapping cone $C_f$, where $f: S\rightarrow X$ is an imbedding of a
sphere of codimension at least 2 in $X$.  If the continuous map
induces monomorphism in homology in the relevant dimensions, then the
freedom of $X$ implies that of $X'$.

\medskip\noindent {\it Proof.}  Let $f: S\hookrightarrow X$ be the
imbedding.  Let $I=[0,\ell]$ with $\ell>>1$ to be determined.  Let
$Cyl_f=X\cup_{f\times\{0\}}(S\times I)$ be the mapping cylinder and
$W=Cyl_f\cup_{id\times\{\ell\}}D$ the mapping cone, where $D$ is a
cell of dimension $\dim(S)+1\leq n-1$.

Let $g$ be a metric on $X$.  Let $h_0=f^*(g)$ be the induced metric on
the sphere $S$.  Let $h_1$ be the metric of a round sphere of
sufficiently large radius so that $h_1\geq h_0$.  We endow the
cylinder $S\times I$ with the metric $(1-x)h_0+x h_1+dx^2$ for $0\leq
x\leq 1$ and $h_1+dx^2$ for $1\leq x\leq \ell.$ Call the resulting
metric $W(g,\ell)$.  Let $p:W(g,\ell)\rightarrow I$ be the map
extending the projection to the second factor $S\times I\rightarrow I$
on the cylinder, while $p(X)=0$ and $p(D)=\ell$.  Let $q\leq n-1$ and
let $z$ be a $q$-cycle in $W$.  By the coarea inequality,
$\vol_q(z)\geq \int_1^l\vol_{q-1}(z\cap p^{-1}(x))dx$.  Hence we find
$x_0\in I$ such that
$$\vol_{q-1}(z\cap p^{-1}(x_0))\leq {1\over \ell-1}\vol_q(z).
\eqno{(6.1)}
$$

Let us show that if $X$ admits a systolically free sequence of metrics
$g_j$ then so does $W$.  Suppose a sequence of cycles $z_j$ in
$W(g_j,\ell)$ satisfies $\vol_q(z_j)=o(\sys_q(g_j))$.  Choosing
$\ell=\ell(j)\geq \sys_q(g_j)$, we obtain
$\lim_{j\rightarrow\infty}{\vol_{q-1}}(z\cap p^{-1}(x_0))=0.$ By the
isoperimetric inequality for small cycles (\cf [31], Sublemma
$3.4.B'$), this $(q-1)$-cycle can be filled by a $q$-disk $B^q$ of
volume which also tends to 0.  Here we must choose $\ell$ big enough
as a function of the metric $g_j$ so that the isoperimetric inequality
would apply.  Let $a= (z\cap d^{-1}([0,x_0])) - B^q$ and $b=z-a$.
Note that $[b]=0$ and so $[a]=[z]\not=0$.  The cycle $a$ lies in the
mapping cylinder which admits a distance-decreasing projection to
$(X,g_j)$, hence $\vol_{q}(a)\geq \sys_{q}(X)$ and so $\vol_{q}(z)\geq
\vol_{q}(a)-\vol_{q}(B)\geq \sys_q(X)-o(1)$.  This shows that the
systoles of $W$ are not significantly diminished compared to those of
$X$.

\medskip Now choose a simplicial structure on $W$.  By the cellular
approximation theorem, a continuous map from $X'$ to $W$ can be
deformed to a simplicial map.  As in [1], we can replace it by a map
which has the following property with respect to suitable
triangulations of $X'$ and $W$: on each simplex of $X'$, it is either
a diffeomorphism onto its image or the collapse onto a wall of
positive codimension.  Let $p$ be the maximal number of $n$-simplices
of $X'$ mapping diffeomorphically to an $n$-simplex of $X\subset W$.
Since the $(n-1)$-cell does not contribute to $n$-dimensional volume,
the pullback of the metric on $W$ is a positive quadratic form on $X'$
whose $n$-volume is at most $p$ times that of $X$.  This form is
piecewise smooth and satisfies natural compatibility conditions along
the common face of each pair of simplices.

\medskip Note that if a smooth compact $n$-manifold $X$ admits
systolically free piecewise smooth metrics, then it also admits
systolically free smooth metrics.  To construct a smooth metric from a
piecewise smooth one, we proceed as in [1].  Given a piecewise smooth
metric $g$, compatible along the common face of each pair of adjacent
simplices, we choose a smooth metric $h$ on $X$ such that $h > g$ at
every point (in the sense of lengths of all tangent vectors).  Let $N$
be a regular neighborhood of small volume of the $(n-1)$-skeleton of
the triangulation.  Choose an open cover of $X$ consisting of $N$ and
the interiors $U_i$ of all $n$-simplices.  Using a partition of unity
subordinate to this cover, we patch together the metrics $g|_{U_i}$
and $h|_N$.  The new metric dominates $g$ for each tangent vector to
$M$.  In particular, the volume of a cycle is not decreased.  
Meanwhile, $n$-dimensional volume is
increased no more than the volume of the regular neighborhood.

\medskip
The piecewise smooth metric on $X'$ may {\it a priori} not be
compatible with its smooth structure, since the triangulation may not
be smooth.  To clarify this point, denote the triangulation by $s$,
and the piecewise smooth metric by $g$.  Consider a smooth
triangulation $s'$, and approximate the identity map of $X'$ by
simplicial map with respect to the two triangulations $s'$ and $s$.
Now we pull the metric $g$ back to $s'$.  This gives a metric $g'$
adapted to the smooth triangulation $s'$, to which we may apply
the argument with the regular neighborhood $N$.

\medskip We have thus obtained a smooth positive form on $X'$.  We
make it definite without significantly increasing its volume by adding
a small multiple of a positive definite form.  The lower bounds for the
systoles are immediate from the injectivity of the map $X'\rightarrow
W$ on the relevant homology groups.

\medskip\noindent {\bf Proposition 6.2.}  Suppose $\C P^2$ and
$S^2\times S^2$ admit 2-free metrics.  Then so does every simply
connected 4-manifold $X$.

\medskip\noindent {\it Proof.}  Let $b=b_2(X)$.  Let $W$ be the
4-skeleton of the Cartesian product $\C P^2\times\ldots\times\C P^2$
($b$ times), where $\C P^2$ comes with its standard 3-cell structure.
Mote that $W$ is the 4-skeleton of the standard model of the
Eilenberg-Maclane space $K(Z^b,2)$.

Choose a CW structure on (the homotopy type of) $X$ whose 2-skeleton
is the wedge of $b$ copies of $S^2$, with a single 4-cell attached.
The identification of the 2-skeleta of $X$ and $W$ extends across the
4-cell since $\pi_3(W)=\pi_3(\C P^2)^b=0$ from the long exact sequence
of the Hopf fibration over $\C P^2.$ The freedom of $W$ and then that
of $X$ is established along the lines of the proof of Lemma 6.1, using
the fact that two distinct 4-cells of $W$ meet along cells of
codimension at least 2.

\medskip\noindent {\bf Proposition 6.3.} Every simply connected
4-manifold $X$ admits a map to a connected sum of copies of $\C P^2$
with either orientation, which induces monomorphism in homology of
dimension 2.

\medskip\noindent {\it Proof.}  Let
$f(\overline{x},\overline{y})=\sum_{i,j=1}^{b_2} f_{i,j}x_iy_j$ be the
intersection form of $X$, where $b_2= b_2(X)$.  Let
$\overline{x},\overline{y} \in H_2(X,Z)$ and let $\sigma=\sigma(X)$ be
the signature of $X$. Let $p=(b_2 + \sigma)/2$ and $q=(b_2 -
\sigma)/2$ and define a manifold $N$ by setting
$$
N=N(b_2,\sigma)=(CP^2\#...\#CP^2)\#(\overline{CP}^2\#...\#\overline{CP}^2),
\eqno{(6.2)}
$$ 
the connected sum of $p$ copies of $CP^2$ and and $q$ copies of
$\overline{CP}^2$.  The intersection form of $N$ is
$$ g(\overline{u},\overline{v})=
u_1v_1+...+u_pv_p-u_{p+1}v_{p+1}-...-u_{b_2}v_{b_2}.
\eqno{(6.3)}
$$ If $pq\not=0$, then $f$ and $g$ are rationally equivalent by the
classification theorem of integer unimodular forms [45].  If, say,
$q=0$, we consider the pair of forms $f+^.\lbrace{-1\rbrace}$ and
$g+^.\lbrace{-1\rbrace}$ and use the `Witt lemma' [16] to the effect
that if two forms are rationally equivalent after adding a common
summand, then they are rationally equivalent.  The rational
equivalence entails the existence of a non-singular integer matrix $C$
such that
$$f(C\overline{u},C\overline{v})=
{\Delta}^2g(\overline{u},\overline{v}),
\eqno{(6.4)}
$$ where $\Delta$ is a nonzero
integer.
It is well known [40] that every simply connected 4-manifold $X$
with intersection form $f=(f_{ij})$ is homotopy equivalent to the complex
$${\vee}^{b_2}_{i=1}S^2_i \hskip3pt{\cup}_{\phi}\hskip3ptD^4,
\eqno{(6.5)}
$$ where
the glueing map $\phi: S^3 \to ({\vee}^{b_2}_{i=1}S^2_i)$ represents
the homotopy class
$${1\over 2}\sum_{i,j=1}^{b_2}f_{i,j}{\lbrack}e_i,e_j{\rbrack}.
\eqno{(6.6)}
$$ Here the $e_i\in {\pi}_2({\vee}^{b_2}_{i=1}S^2_i)$ are the
canonical generators, while the Whitehead products $[e_i,e_j]$ for
$i<j$ together with ${1\over 2}[e_i,e_i]$ form a basis for
$\pi_3({\vee}^{b_2}_{i=1}S^2_i)$ by the Hilton-Milnor theorem (\cf
[37], [47]).

Let us consider the map
$${\xi}_C: {\vee}^{b_2}_{i=1}S^2_i \rightarrow
{\vee}^{b_2}_{i=1}S^2_i
\eqno{(6.7)}
$$ defined by the matrix $C$ introduced above.
The usual calculation shows that
$$\matrix{ {{\xi}_C}_*(\lbrace{\phi\rbrace}) & = &
{{\xi}_C}_*(1/2\sum_{i,j=1}^{b_2}f_{i,j} {\lbrack}e_i,e_j{\rbrack})
\hfill\cr &=& 1/2\sum_{i,j=1}^{b_2} f_{i,j}{\lbrack}
{{\xi}_C}_*(e_i),{{\xi}_C}_*(e_j){\rbrack} \hfill\cr &=&
{\Delta}^2/2({\lbrack}\overline{e_1},\overline{e_1}{\rbrack}+\ldots+{\lbrack}
\overline{e_p},\overline{e_p}{\rbrack}-
{\lbrack}\overline{e_{p+1}},\overline{e_{p+1}}{\rbrack}-\ldots-{\lbrack}
\overline{e_{b_2}},\overline{e_{b_2}}{\rbrack}), \hfill\cr}
\eqno{(6.8)}$$ where the $\overline{e_{i}} \in {\pi}_2(N)$ are the
canonical generators defined by the inclusion
$${\vee}_{i=1}^{b_2} S^2_i \subset N.
\eqno{(6.9)}
$$
The computation shows that ${\xi}_C$ can be extended to a map
$${\Xi}_C:X =_{hom} {\vee}^{b_2}_{i=1}S^2_i
\hskip3pt{\cup}_{\phi}\hskip3ptD^4\rightarrow N , \eqno(6.10)$$
moreover $deg\hskip2pt{\Xi}_C={\Delta}^2$.  

\medskip\noindent Note that Proposition 6.3 provides an alternative
proof, not involving CW complexes, of the fact that if $\C P^2$ is
free then every simply connected 4-manifold $X$ is free.  Here we
apply Lemma 6.1 to the mapping cone of the empty map.

\medskip\noindent {\bf Proposition 6.4.} The following three
assertions are equivalent:

(i) $\C P^2$ is 2-free;

(ii) $\C P^2\#\cb$ is 2-free;

(iii) $S^2\times S^2$ is 2-free.
\medskip\noindent
{\it Proof.} 
$(i)\Rightarrow (ii)$: If $\C P^2$ has free metrics then so does
$\cb$.  Taking a connected sum by a thin long tube
produces a metric on  $\C P^2\#\cb$ which admits
distance-decreasing projections of degree 1 to each of the summands.
Since 
$H_2(\C P^2\#\cb)
=H_2(\C P^2)+ H_2(\cb)=\Z+ \Z$,
the implication follows.
\medskip\noindent
$(ii)\Rightarrow (iii)$: Recall that 
$\C P^2\#\cb$ is the nontrivial 2-sphere bundle
over $S^2$.  Let $f:S^2\rightarrow S^2$ be a degree 2 map.
The pullback of the nontrivial bundle by $f$
is the trivial one (\cf the $w_2$ discussion preceding (5.1)).
We thus obtain a map 
$f_*:S^2\times S^2\rightarrow \C P^2\#\cb$
inducing an injective homomorphism
$H_2(S^2\times S^2)\rightarrow H_2(\C P^2\#\cb)$.
We approximate this map by a simplicial one as in [1],
and pull back the free metrics from
$\C P^2\#\cb$ to $S^2\times S^2$ (\cf  6.5),
proving this implication.
\medskip\noindent
$(iii)\Rightarrow (i)$:
We cannot map $\C P^2$ to $S^2\times S^2$ in such a way as to induce
a monomorphism in $H_2$.
The obstruction lies in the group $\pi_3(S^2\times S^2)=\Z+\Z$.
We can, however, kill two birds with one stone, or more precisely
with one 3-cell, eliminating 3-dimensional homotopy (at least
rationally).

We glue in a 3-ball $B^3$ to $S^2\times S^2$ along 
the diagonal sphere to obtain a CW complex $W$
which admits a map from $\C P^2$ inducing an
injective homomorphism on 2-di\-mensional homology.
The 2-systole of $W$ obeys the same asymptotics
as that of the area-rich metrics on $S^2\times S^2$.
The coveted metrics on $\C P^2$ are pulled back from 
$W$ by the map $\C P^2\rightarrow W$.
\medskip\noindent
Let $W=(S^2\times S^2)\cup B^3$ where
the 3-ball is glued in along the imbedded diagonal
2-sphere representing the element
$(1,1)\in \pi_2(S^2\times S^2)=\Z+\Z.$
Let $g: S^2\times S^2\rightarrow W$ be the inclusion.
We will specify a map 
$f:\C P^1\rightarrow S^2\times S^2$ such that the induced 
homomorphism $(g\circ f)_2:\pi_2(\C P^1)\rightarrow\pi_2(W)$
is injective while
$(g\circ f)_3:\pi_3(\C P^1)\rightarrow\pi_3(W)$ is zero.
Recall that $\C P^2=\C P^1\cup_h B^4$,
where $h$ is the generator of $\pi_3(\C P^1)$
(in fact, $W$ is homotopy equivalent to
$S^2\cup_{2h} B^4$, whence the existence of the map
extending a degree 2 map on the 2-sphere).
Therefore the map $g\circ f$ extends to a map
$\C P^2\rightarrow W$ which induces an injective homomorphism
in $\pi_2$, and hence in $H_2$ by the Hurewicz theorem
(\cf [25]).
The property $(g \circ f)_3=0$ follows from the fact that
$\Image(g_3)$ is 2-torsion while $\Image(f_3)$ is even
(see Lemma 6.6 below).
The proof is completed by applying Lemma 6.1 to the mapping cone
of the inclusion of the diagonal in $S^2\times S^2$.

\medskip\noindent
{\bf Lemma 6.6.}  Let $f:\C P^1\rightarrow S^2\times S^2$
be a map sending $\C P^1$ to the first factor with degree 2.
Let $W=(S^2\times S^2)\cup B^3$ where
the 3-ball is glued in along the imbedded diagonal
2-sphere representing the element
$(1,1)\in \pi_2(S^2\times S^2)=\Z+\Z.$
Let $g: S^2\times S^2\rightarrow W$ be the inclusion.
Then the map $g\circ f$ extends to a map
$\C P^2\rightarrow W$.
\medskip\noindent
{\it Proof.}
Let $e_1, e_2$ be the standard generators of 
$\pi_2(S^2\times S^2)=\Z+\Z $.  The group
$\pi_3(S^2\times S^2)=\Z+\Z$ is generated  by
elements $h_1$ and $h_2$ satisfying
$2h_i=[e_i,e_i]$ for $i=1,2$, where the brackets
denote the Whitehead product
(\cf [25], p.\ 74 and [47], p.\ 495, theorem 2.5).
Here $[e_1,e_2]=0$ in $\pi_3(S^2\times S^2)$ by definition
of the Whitehead product.
We have $g_2(e_1+e_2)=0\in \pi_2(W)$ by definition of $W$.

Now let $e\in\pi_2(\C P^1)$ and $h\in\pi_3(\C P^1)$ be the respective
generators, so that $[e,e]=2h$.  Consider the class
$2e_1\in\pi_2(S^2\times S^2)$, and take a representative $f\in 2e_1$.
Then $f_2(e)=2e_1\in\pi_2(S^2\times S^2)$.  By naturality of the
Whitehead product,
$$2f_3(h)=[f_2(e),f_2(e)]=[2e_1,2e_1]=4[e_1,e_1]=8h_1
\in\pi_3(S^2\times S^2)=\Z+\Z.  \eqno{(6.12)} $$ In a free abelian
group we can divide by 2, obtaining $$f_3(h)=4h_1\in\pi_3(S^2\times
S^2).  \eqno{(6.13)} $$ Now the lemma follows from the fact that $\C
P^2=\C P^1\cup_h B^4$ and the following calculation:
$$\matrix{ (g\circ f)_3(h) & = & g_3(4h_1)=2g_3(2h_1)=
2[g_2(e_1),g_2(e_1)] \hfill\cr &=& 2[g_2(e_1+e_2),g_2(e_1)]=
2[0,g_2(e_1)]=0.  \hfill\cr} \eqno{(6.14)}$$ 

\bigskip\noindent
\vfill\eject 
\noindent {\large References} 

\medskip\noindent [1]
I. Babenko, Asymptotic invariants of smooth manifolds, {\it Russian
Acad.\ Sci.\ Izv.\ Math.\/} 41 (1993) 1-38.  

\medskip\noindent [2] I. Babenko, M. Katz, and A. Suciu, Volumes,
middle-dimensional systoles, and Whitehead products, preprint (1997).

\medskip\noindent [3] V. Bangert, Minimal geodesics, {\it Ergod.\ Th.\
and Dynam.\ Sys.\/} 10 (1990) 263-286.

\medskip\noindent [4] L. B\'erard Bergery \& M. Katz, Intersystolic
inequalities in dimension 3, {\it Geometric and Functional Analysis}
4, No.\ 6 (1994) 621-632.

\medskip\noindent [5] L. B\'erard Bergery \& M. Katz, Homogeneous
freedom in middle dimension, in progress.

\medskip\noindent [6] M. Berger, Du c\^ot\'e de chez Pu, {\it Ann.\
Scient.\ Ec.\ Norm.\ Sup.\/} 5 (1972) 1-44.

\medskip\noindent [7] M. Berger, A l'ombre de Loewner {\it Ann.\
Scient.\ Ec.\ Norm.\ Sup.\/} 5 (1972) 241-260.

\medskip\noindent [8] M. Berger, Systoles et applications selon
Gromov, expos\'e 771, S\'eminaire N. Bourbaki, Ast\'erisque 216 (pp.\
279-310), Soci\'et\'e Math\'ematique de France, 1993.

\medskip\noindent [9] M. Berger, private communication (1 March
1994).

\medskip\noindent [10] M.  Berger, Riemannian geometry during the
second half of the twentieth century, I.H.E.S. preprint (March 1997).

\medskip\noindent [11] A. Besicowitch, On two
problems of Loewner, {\it J. London Math.\ Soc.\/} 27 (1952) 141-144.

\medskip\noindent [12] R. Bott, The stable homotopy of the
classical groups, {\it Ann.\ Math.\/} 70 (1959) 313-337.

\medskip\noindent [13] Yu. Burago and V. Zalgaller, {\it Geometric
inequalities.\/} Springer, 1988.

\medskip\noindent [14] M. do Carmo, {\it Riemannian geometry,}\/
Birkhauser, 1992.

\medskip\noindent [15] E. Cartan, {\it La topologie des espaces
representatifs des groupes de Lie.} Act.\ Sci.\ Ind., n.\ 358,
Hermann, Paris 1936.

\medskip\noindent [16] J. Cassels, {\it Rational quadratic forms.\/}
Academic Press, 1978.

\medskip\noindent [17] B.~Colbois and J. Dodziuk, Riemannian Metrics
with Large $\lambda_1$, {\it Proc.\ Amer.\ Math.\ Soc.,}\/ {\bf 122}
(1994) 905--906.

\medskip\noindent [18] G. Courtois, Exemple de vari\'et\'e
riemannienne et de spectre, in Rencontres de th\'eorie sp\'ectrale et
g\'eom\'etrie, Grenoble 1991.

\medskip\noindent [19] J. Diedonne, {\it A history of algebraic and
differential topology 1900-1960.\/} Birkhauser, 1989.

\medskip\noindent [20] J. Dodziuk, Nonexistence of universal upper
bounds for the first positive eigenvalue of the Laplace-Beltrami
operator, Proceedings of the 1993 Joint Summer Research Conference on
Spectral Geometry, {\it Contemporary Mathematics,}\/ {\bf 173} (1994)
109-114.

\medskip\noindent [21] B. Dubrovin, A. Fomenko, S. Novikov, {\it
Modern geometry--methods an applications. Part II. The geometry and
topology of manifolds.} Graduate Texts in Mathematics 104.  Springer,
1985.  

\medskip\noindent [22] Encyclopedic dictionary of mathematics, Second
Edition (Japanese encyclopedia).  The MIT Press, 1987.

\medskip\noindent [23] H. Federer, Real flat chains, cochains and
variational problems, {\it Indiana Univ.\ Math.\ J.\/} 24 (1974)
351-407.

\medskip\noindent [24] H. Federer \& W. Fleming, Normal and integral
currents, {\it Ann.\ of Math.\/} 72 (1960) 458-520.

\medskip\noindent [25] A. Fomenko, D. Fuchs, and V. Gutenmacher.
{\it Homotopic topology.\/} Akad\'{e}miai Kiad\'{o}, Budapest, 1986.

\medskip\noindent [26] S. Gallot, D. Hulin, \& J. Lafontaine, {\it
Riemannian Geometry\/}, Springer, 1987.

\medskip\noindent [27] M.  Glezerman and L. Pontryagin, {\it
Intersections in manifolds,\/} Amer Math. Soc. Translations No. 50,
1951.  (Translation of Uspehi Matematicheskih Nauk (N.S.) 2, no. 1
(17) (1947) 50-155.)

\medskip\noindent [28] C. Godbillon, {\it Topologie alg\'ebrique.\/}
Hermann, Paris, 1971.  

\medskip\noindent [29] M. Greenberg \& J. Harper, Algebraic topology:
a first course. Addison-Wesley, 1981.

\medskip\noindent [30] M. Gromov, {\it Structures m\'etriques pour les
vari\'et\'es riemanniennes\/} (edited by J. Lafontaine, and P. Pansu).
Cedic, 1981.

\medskip\noindent [31] M. Gromov, Filling Riemannian manifolds, {\it
J. Diff.\ Geom.\/} 18 (1983) 1-147.

\medskip\noindent [32] M. Gromov, Systoles and intersystolic
inequalities, I.H.E.S.\ preprint (1992).  

\medskip\noindent [33] M.  Gromov, Metric invariants of Kahler
manifolds, in {\it Proc.\ Differential Geometry and Topology,}\/
R. Caddeo, F. Tricerri ed., World Sci., 1993.

\medskip\noindent [34] M. Gromov, Systoles and intersystolic
inequalities, pp.\ 291-362, in A. Besse, Actes de la table ronde de
g\'eom\'etrie diff\'erentielle en l'honneur de Marcel Berger
(S\'eminaires et Congr\`es 1, Soci\'et\'e Math\'ematique de France,
1996).

\medskip\noindent [35] M. Gromov, {\it The metric geometry of
Riemannian and non-\-Riemannian spaces} (edited by J. Lafontaine and
P. Pansu).  Birkhauser, 1997 (expanded English edition of the book
[30]).

\medskip\noindent [36] J. Hebda, The collars of a Riemannian
manifold and stable isosystolic inequalities, {\it Pacific J. Math.\/}
121 (1986) 339-356.

\medskip\noindent [37] P. Hilton, On the homotopy groups of the
union of spheres, {\it J.  London Math.\ Soc.\/} 30 (1955) 154-172.

\medskip\noindent [38] M. Katz, Counterexamples to isosystolic
inequalities, {\it Geometriae Dedicata\/} 57 (1995) 195-206.

\medskip\noindent [39] S. Lefschetz, Intersections and transformations
of complexes and manifolds, {\it Trans.\ Amer.\ Math.\ Soc.\/} 28
(1926) 1-49.

\medskip\noindent [40] R. Mandelbaum, {\it Four-dimensional
Topology.\/}

\medskip\noindent [41] J. Milnor and J. Stasheff,
Characteristic classes.  Princeton University Press, 1974.

\medskip\noindent [42] C. Pittet, Systoles on $S^1\times S^n$, {\it
Diff.\ Geom.\ Appl.\/} 7 (1997) 139-142.

\medskip\noindent [43] P. Pu, Some inequalities in certain
nonorientable manifolds, {\it Pacific J. Math.\/} 2 (1952) 55-71.

\medskip\noindent [44] J. Schwarz, {\it Differential geometry and
Topology.\/} Gordon and Breach, New York, 1968.

\medskip\noindent [45] J.-P. Serre, {\it Cours d'arithm\'etique.\/}
Presses Universitaires de France, 1970.

\medskip\noindent [46] R.  Thom, Quelques propri\'et\'es globales des
vari\'et\'es diff\'erentiables, {\it Comment.\ Math.\ Helv.\/} 28
(1954) 17-86.

\medskip\noindent [47] G. Whitehead, {\it Elements of homotopy
theory.\/} Springer, 1978.

\bigskip \vfill \noindent addresses:

\medskip\noindent D\'epartement des Sciences Math\'ematiques,
Universit\'e de Montpellier 2, place Eug\`ene Bataillon, 34095
Montpellier, France.  \ \ email: babenko@darboux.math.univ-montp2.fr

\medskip\noindent UMR 9973, D\'epartement de Math\'ematiques,
Universit\'e de Nancy 1, B. P. 239, 54506 Vandoeuvre, FRANCE.\ \
email: katz@iecn.u-nancy.fr

\vfill\eject\end